\documentclass[11pt]{amsart}

\usepackage[letterpaper,top=1in, left=1in, right=1in, bottom=1in]{geometry}

\usepackage{todonotes}
\usepackage{xcolor}

\usepackage{amssymb,amsmath,epsfig,color}

\usepackage{subcaption}
\usepackage[abs]{overpic}

\usepackage[hidelinks,pagebackref,pdftex]{hyperref}
\usepackage{xcolor}

\newtheorem{theorem}{Theorem}[section]

\newtheorem{proposition}[theorem]{Proposition}
\newtheorem{corollary}[theorem]{Corollary}
\newtheorem{conjecture}[theorem]{Conjecture}
\theoremstyle{definition}
\newtheorem{definition}[theorem]{Definition}
\newtheorem{example}[theorem]{Example}
\newtheorem{remark}[theorem]{Remark}

\title{On the notion of Khovanov A-adequacy}

\author[Buchanan]{Lizzie Buchanan}
\address{Department of Mathematics, University of Iowa, USA}
\email{{\rm elizabeth-buchanan@uiowa.edu}}

\author[Guo]{Huizheng Guo}
\address{Department of Mathematics, The George Washington University, Washington DC, USA}
\email{{\rm hguo30@gwu.edu}}

\author[Montoya-Vega]{Gabriel Montoya-Vega}
\address{Department of Mathematics, The Graduate Center CUNY, NY, USA and \newline \indent Department of Mathematics, University of Puerto Rico-R\'io Piedras, San Juan, PR}
\email{{\rm gabrielmontoyavega@gmail.com}}

\author[Rong]{Yongwu Rong}
\address{Department of Mathematics, Queens College and the Graduate Center, CUNY, NY, USA}
\email{{\rm yrong@qc.cuny.edu}}

\author[Silvero]{Marithania Silvero}
\address{Departamento de Álgebra, Universidad de Sevilla, Seville, Spain}
\email{{\rm marithania@us.es}}

\begin{document}

\begin{abstract}
The concept of adequate links, introduced by Lickorish and Thistlethwaite as a generalization of alternating links, has recently gained interest among knot theorists in the context of Khovanov homology. Przytycki and Silvero introduced the more general concept of Khovanov adequacy: a diagram is Khovanov-adequate if its associated Khovanov chain complexes at both potential maximal and minimal quantum gradings have non-trivial homology. This article explores Khovanov adequacy within the framework of independence complexes and the calculation of the homotopy type of extreme Khovanov spectra.
\end{abstract}

\maketitle


\tableofcontents



\section{Introduction}

Adequate links were introduced by W. B. R. Lickorish and M. B. Thistlethwaite as a generalization of alternating links in the context of the proof of the First Tait Conjecture \cite{LT}. As it is typical in knot theory, the adequacy character of a link is defined in terms of the existence of a diagram having certain properties. Namely, a link is adequate if it can be represented by an adequate diagram (see Definition \ref{Defadequate}). 

\ 

A key property of adequate diagrams is the fact that they realize both the {\it{(potential)}} maximal and minimal degrees in their Jones polynomial (i.e., the span of their Jones polynomial is {\it{as big as possible}}). This property was crucial to establish a relation between the number of crossings of a diagram $D$ and the span of the Jones polynomial of the link represented by $D$, which led to the proof of First Tait Conjecture on minimality of crossing number of reduced alternating diagrams. 

\ 

Khovanov homology is a powerful link invariant introduced by M. Khovanov as a categorification of the Jones polynomial \cite{Kho1}. This invariant has a richer algebraic structure and, surprisingly, detects the unknot \cite{KrMr}. Given a link diagram $D$ representing a link $L$, Khovanov built a bigraded chain complex associated to $D$, whose bigraded homology groups are link invariants. These groups $Kh^{i,j}(L)$ are known as Khovanov homology groups of $L$ and the indexes $i$ and $j$ as {\it{homological}} and {\it{quantum}} gradings, respectively. In \cite{PS2} the idea of adequacy was generalized to the concept of Khovanov-adequacy. More precisely, a diagram is Khovanov-adequate if its associated Khovanov chain complexes at both {\it{(potential)}} maximal and minimal quantum gradings give rise to non-trivial Khovanov homology (i.e., the span on the quantum grading of their Khovanov homology is {\it{as big as possible}}). A link is Khovanov-adequate if it can be represented by a Khovanov adequate diagram.

\ 

Since the Jones polynomial of a link can be obtained as the Euler characteristic of its Khovanov homology, it is clear that adequate links are Khovanov-adequate. In fact, this is a proper inclusion, as can be easily shown by using the fact that the Jones polynomial of any adequate link is monic (i.e., the extreme coefficients are equal to $ \pm 1$) and its extreme Khovanov homology (in both, maximal and minimal quantum gradings) is isomorphic to $\mathbb{Z}$. Although it is an open question which groups can appear in the extreme homology of Khovanov-adequate links, there are several results suggesting that torsion groups cannot occur (see \cite{PS1, PS2, OSYY}). 

\ 


In this paper we delve into the understanding of Khovanov adequate links. To do so, we use a characterization based on the geometric realization of extreme Khovanov homology of a link diagram in terms of the independence complex of a certain graph constructed from the diagram \cite{GMS}. This construction  has been proven \cite{CS} to be stably homotopy equivalent to the Khovanov spectrum introduced in \cite{LS} when fixing extreme quantum grading. 

\

The structure of the paper is as follows: In Section \ref{sectcirclegraphs} we recall some preliminaries on circle graphs and their independence complexes and pose the conjecture about their homotopy type. Section \ref{sectKhovade} is devoted to introduce Khovanov adequate diagrams and to state the result relating extreme Khovanov homology of a diagram $D$ with the cohomology of certain independence complex associated to $D$. In Section \ref{Bipartite} we explore independence complexes arising from bipartite circle graphs and provide an obstruction for such independence complexes. Finally, we show some examples in Section \ref{Examples} to illustrate some of the ideas discussed in the paper.

\section{Circle graphs and independence complexes}\label{sectcirclegraphs}

In this section we review the definition of circle graphs and state a conjecture involving the homotopy type of their independence complexes.

\begin{definition}
A chord diagram $\mathcal{C}$ is a circle together with a finite set of chords with disjoint boundary points. The circle graph $G_\mathcal{C}$ associated to the chord diagram $\mathcal{C}$ is the intersection graph of its chords, that is, the simple graph constructed by associating a vertex to each chord in $\mathcal{C}$ and adding an edge connecting two vertices if the corresponding chords intersect.
Given a graph $G$, its independence complex $I(G)$ is the simplicial complex whose simplices are the subsets of pairwise non-adjacent vertices of $G$.
\end{definition}

\begin{figure}[h]
    \centering
\includegraphics[width=0.6\linewidth]{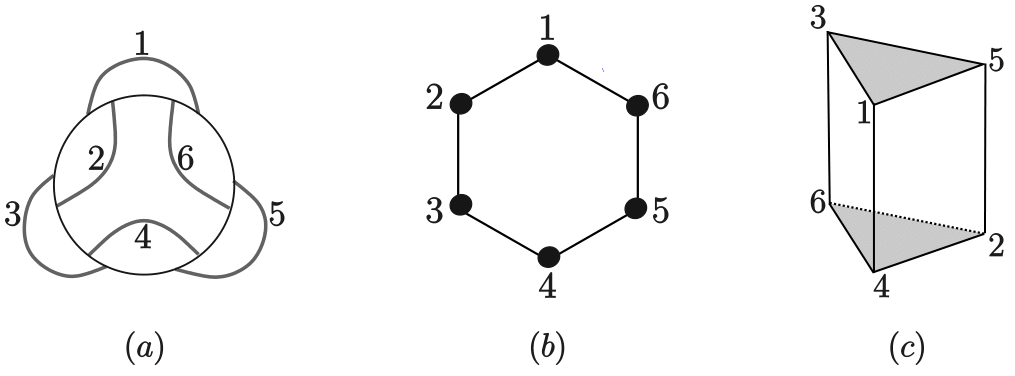}
\caption{A chord diagram $\mathcal{C}$, its associated circle graph $G_\mathcal{C}$ and the geometric realization of the corresponding independence complex $I(G_\mathcal{C}$) are shown in (a), (b) and (c), respectively.}
\label{Examplecircle}
\end{figure}

It is not hard to see that not all graphs are circle graphs, that is, some graphs cannot be represented as the intersection graphs of any chord diagram (see \cite{Bou} for a theoretical characterization of circle graphs). In this paper we are particularly interested in those circle graphs which are bipartite, the reason will be clear in the next section. 

\ 

It is natural to ask which complexes, up to homotopy type, can arise as the independence complex of a circle graph. In \cite[Proposition~4.3]{PS1} it was given an explicit method to construct a circle graph $G$ whose independence complex $I(G)$ has the homotopy type of any given wedge of spheres. The following conjecture concerns the other implication. 

\begin{conjecture}\cite{PS1}\label{conj1}
The independence complex associated to a circle graph is homotopy equivalent to a wedge of spheres. 
\end{conjecture}

The conjecture has been proved to be true for several families of graphs such as paths, trees, polygons, cactus graphs, outerplanar graphs, permutation graphs, non-nested graphs, and those circle graphs arising as the Lando graphs of several families of links (see \cite{Koz, PS1, PS2, OSYY}, for example). We finish this section with some observations and results that will be used throughtout the paper:

\begin{enumerate}
    \item If $G$ is the empty graph, we set $I(G)$ to be $S^{-1}$, that we call \textit{sphere of dimension -1} and whose suspension is $S^0$ (i.e., $\Sigma S^{-1} = S^0$).
    \item Given two chord diagrams $\mathcal{C}_1$ and $\mathcal{C}_2$, write $G_1$ and $G_2$ for their associated circle graphs, respectively. Then, the circle graph $G_{\mathcal{C}_1 \sqcup \mathcal{C}_2} = G_1 \sqcup G_2$, and therefore $I(G_{\mathcal{C}_1 \sqcup \mathcal{C}_2}) = I(G_1) \ast I(G_2)$, where $\ast$ denotes the join operation.
\end{enumerate}

\begin{proposition}\cite{Koz}\label{cycle pattern lemma} Let $L_n$ be the path on $n+1$ vertices and $C_n$ the polygon on $n$ edges. Then
$$I(L_n) \sim_h \left\{ \begin{array}{lll} \bullet & & \mbox{ if } n=3k,  \\ \\ S^k & & \mbox{ if } n=3k+1, 3k+2; \end{array} \right. \quad \quad I(C_n) \sim_h \left\{ \begin{array}{lcc} S^{k-1} & & \mbox{ if } n=3k \pm 1,  \\ \\ S^{k-1} \vee S^{k-1} & & \mbox{ if } n=3k. \end{array} \right.$$

\end{proposition}

\vspace{0.1cm}

\section{Khovanov adequate diagrams}\label{sectKhovade}

\begin{figure}[h]
    \centering
\includegraphics[width=0.7\linewidth]{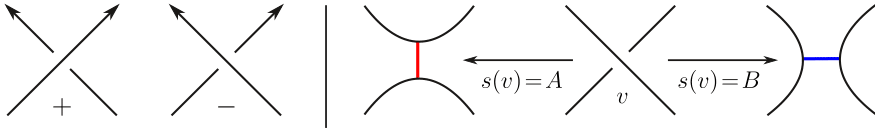}
\caption{Left: a positive (+) and a negative (-) crossing. Right: the smoothing of a crossing, according to its $A$ or $B$ label.}
\label{Cruces}
\end{figure}

Let $D$ be a link diagram with $p$ positive and $n$ negative crossings (see Figure \ref{Cruces} for sign convention). A \textit{Kauffman state} $s$ is an assignation of a label, $A$ or $B$, to each crossing of $D$, that is, $s: cr(D) \to \{A,B\}$. We write $sD$ for the set of (topological) circles and (colored) chords obtained by smoothing each crossing of $D$ according to its label following Figure \ref{Cruces}. Denote by $s_A$ (resp. $s_B$) the state assigning an $A$-label (resp. a $B$-label) to every crossing of $D$. 

\begin{definition}\label{Defadequate}
A diagram $D$ is said to be $A$-adequate (resp. $B$-adequate) if every chord in $s_AD$ (resp. in $s_BD$) has both endpoints in different circles. We say that $D$ is adequate if it is both $A$-adequate and $B$-adequate. A link is said to be ($A$/$B$-)adequate if it admits an ($A$/$B$-)adequate diagram. 
\end{definition}

An enhancement of an arbitrary state $s$ is a map $x$ which assigns a sign $\pm 1$ to each of the circles in $sD$. Write $\tau(s,x) = \sum x(c)$, where $c$ ranges over all circles in $sD$, and define, for the enhanced state $(s,x)$, the integers $$i(s,x) = i(s) = -n + |s^{-1}(B)| \quad \quad \mbox{and} \quad \quad j(s,x) = p - 2n + |s^{-1}(B)| + \tau(s,x).$$

The set of enhanced states associated to $D$ are the generators of the chain groups in the Khovanov complex of $D$. More precisely, the complex has the form $$\ldots \, \longrightarrow \, C^{i,j}(D) \, \stackrel{\partial_i}{\longrightarrow} \, C^{i+1,j}(D) \, \longrightarrow \ldots,$$
where $C^{i,j}(D)$ is the free abelian group generated by the set of enhanced states $(s,x)$ of $D$ with $i(s) = i$ and $j(s,x) = j$ and $\partial_i$ is Khovanov's differential (see \cite{Vir} for a precise definition). Khovanov proved that the associated homology groups $Kh^{i,j}(D)$ are link invariants, and they are referred to as {\it{Khovanov homology groups}} of the link represented by $D$.

\ 

The Khovanov complex is not a link invariant and it is defined for every $i,j \in \mathbb{Z}$. However, for a fixed diagram $D$ there is an {\it{extreme}} value $j_{\min}(D) \in \mathbb{Z}$ such that all chain groups $C^{i,j}(D)$ become trivial when $j<j_{\min}(D)$. Let $$j_{\min}(D) = \min \{j(s,x) \, | \, (s,x) \mbox{ is an enhanced state of }D\},$$ and denote the complex $\{C^{*,j_{\min}}(D), \partial_*\}$ the \textit{extreme Khovanov complex} and its associated homology groups $Kh^{*, j_{\min}}(D)$ the \textit{extreme Khovanov homology groups} of $D$. We write $$\underline{j}(L) = \min \{j \, | \, Kh^{*,j}(L) \neq 0\}.$$ 

Observe that $j_{\min}(D) \leq \underline{j}(L)$ for every diagram $D$ representing $L$.  

\begin{definition}\label{DefKhovadequate}
A diagram $D$ representing a link $L$ is Khovanov $A$-adequate if $j_{\min}(D) = \underline{j}(L)$, that is, if $Kh^{*,j_{\min}(D)}(L)$ is non-trivial.
\end{definition}

It is clear from the definition of quantum grading that $j_{\min}(D) = j(s_A,x_-)$, where $x_-$ the constant map associating a negative sign to every circle in $s_AD$. Therefore, $j_{\min}(D) = p-2n-|s_AD|$, where $|s_AD|$ denotes the number of circles in $s_AD$.

\ 

In \cite{GMS} authors introduced a geometric realization of extreme Khovanov homology in terms of certain simplicial complex, which was later shown to be stably homotopy equivalent to the celebrated Khovanov spectrum introduced by Lipshitz and Sarkar in \cite{LS},  when restricted to extreme quantum grading (see \cite{CS}). We review the construction now.

\begin{definition}
Given a link diagram $D$, its Lando graph $G_D$ is constructed from $s_AD$ by associating a vertex to each chord having both endpoints in the same circle, and by adding an edge connecting two vertices if the endpoints of the corresponding chords alternate along the same circle. 
\end{definition}

We can think of the Lando graph $G_D$ as the circle graph associated to the chord diagram obtained from $s_AD$ after removing those chords having their endpoints in different circles. In fact, the family of Lando graphs coincides with the family of bipartite circle graphs. 

\begin{theorem}\cite{GMS}
Let $D$ be a link diagram with $n$ negative crossings, and let $I(D)$ be the independence complex associated to its Lando graph $G_D$. Then, the simplicial cochain complex of $I(D)$ is isomorphic to the extreme Khovanov complex $\{C^{i,j_{\min}}(D), \partial_i\}$, shifted by $n-1$. As a consequence, $$Kh^{i,j_{\min}}(D) \simeq H^{i-1+n}(I(D)).$$   
\end{theorem}

The above theorem allows us to redefine Khovanov $A$-adequate diagrams as those $D$ such that $I(D)$ is not contractible. If Conjecture \ref{conj1} were true, it would imply that extreme Khovanov homology of any diagram is torsion-free, leading to an obstruction for a link to be Khovanov $A$-adequate. In fact, since we are just interested in circle graphs arising as Lando graphs associated to any diagram, we can restrict our study to the family of bipartite graphs, and consider the following conjecture instead.

\begin{conjecture} \cite{PS1}\label{conj2}
The independence complex associated to a bipartite circle graph is homotopy equivalent to a wedge of spheres.
\end{conjecture}

In the following section we explore independence complexes arising from bipartite circle graphs. 

\section{Bipartite graphs and their independence complexes}
\label{Bipartite}

This section is motivated by Conjectures \ref{conj1} and \ref{conj2} stated in the previous sections. Recall that, in the setting of extreme Khovanov homology, the corresponding independence complexes arise from a graph that is both bipartite and a circle graph. 
Independence complexes arisen from bipartite circle graphs were studied by Nagel and Reiner \cite{NR} and by Jonsson \cite{Jon}. In particular, they characterized when a simplicial complex is the independence complex of a bipartite circle graph.  The following is a version of the result in \cite{Jon}.

\begin{theorem}[\cite{Jon} Theorem 1]\label{Jonsson}
\label{JonssonThm}
A simplicial complex $X$ is homotopy equivalent to the independence complex of some bipartite graph $G$ if and only if $X$ is homotopy equivalent to the suspension of some simplicial complex.   
\end{theorem}

As a quick remark, we note that the suspension space is always connected except for the suspension of the empty set (that we called $S^{-1}$), which is $S^0$. Thus we get:

\begin{corollary}
The independence complex of a graph is either connected or it has the homotopy type of $S^0$. 
\end{corollary}

Given a simplicial complex $Y$, Jonsson gave an explicit method to obtain a graph $G$ so that $I(G) \sim_h \Sigma Y$. We describe his construction: Denote by $V(Y)$ the vertex set of $Y$ and by $M(Y)$ the set of maximal faces of $Y$. We construct a bipartite graph $G = G(Y)$ as follows: the vertex set of $G$ is the disjoint union $V(Y) \sqcup M(Y)$, and two vertices $v \in V(Y)$ and $w \in M(Y)$ are connected by an edge in $G$ if and only if $v \notin w$. Notice that we use the same symbols for the vertices in $G$ as for the simplices in $Y$, with a slight abuse of notation.

\begin{figure}[ht]
\centering
\begin{subfigure}{.49\textwidth}
\centering		
\includegraphics[width=0.45\linewidth]{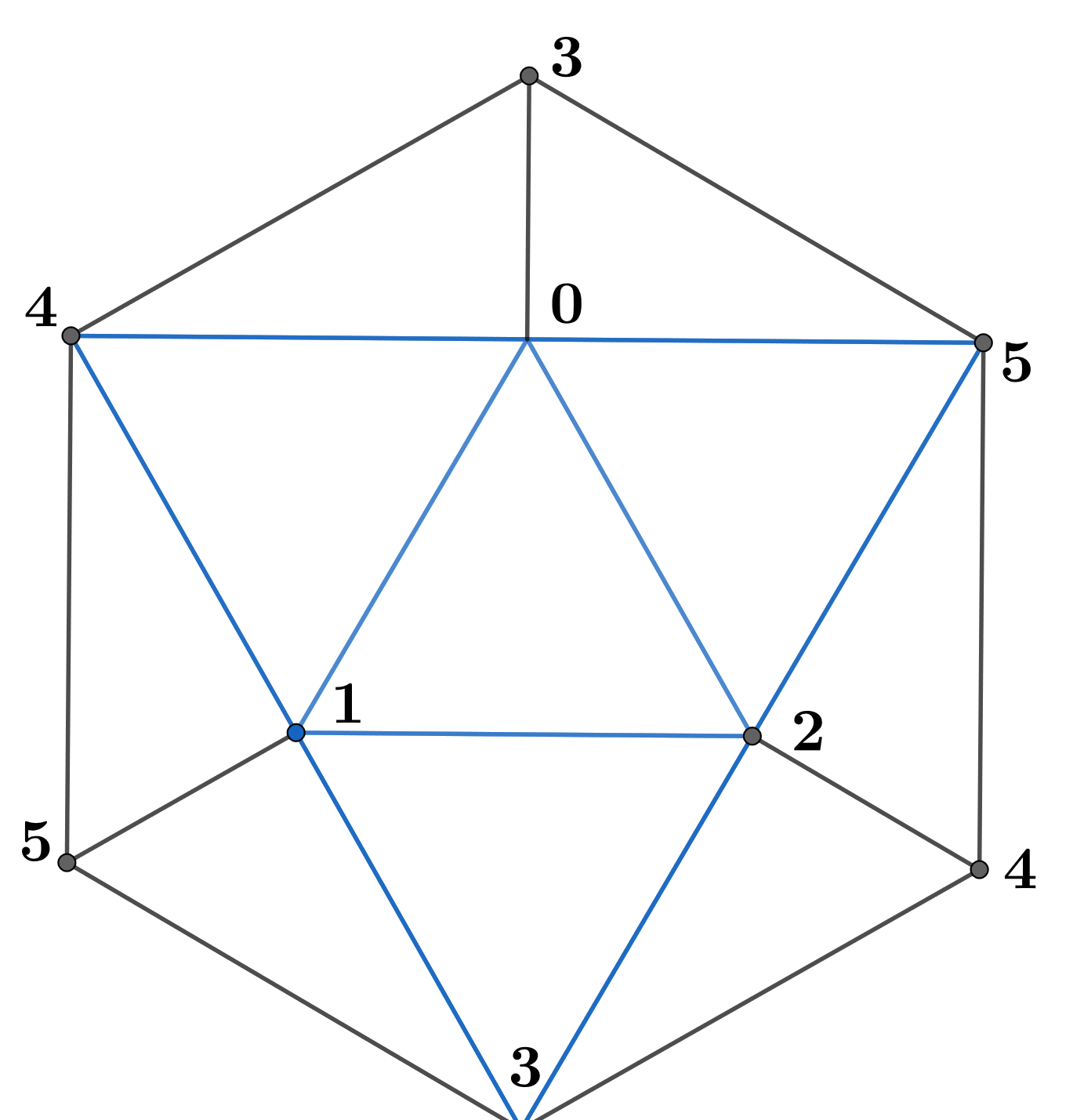}
\caption{A triangulation of $\mathbb{R}P^2$.}
\label{fig:RP2Triangulation}
\end{subfigure}
\begin{subfigure}{.49\textwidth}
\centering
\includegraphics[width=0.5\linewidth]{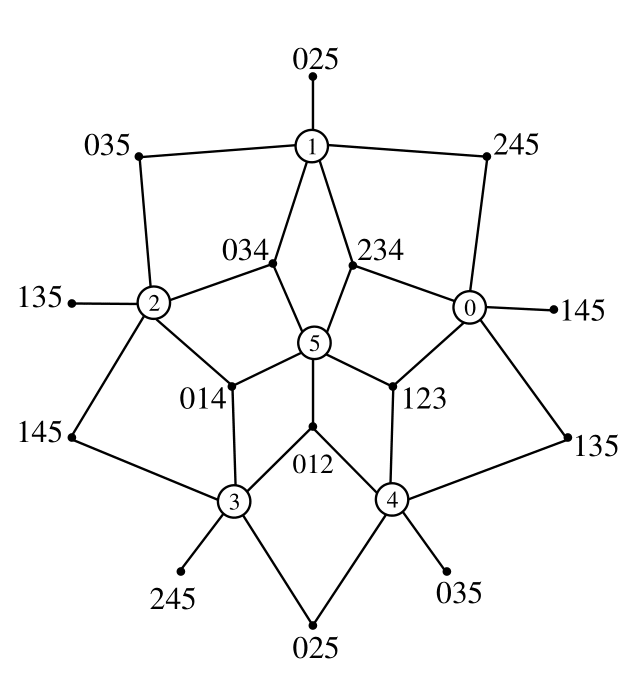}
\caption{The associated bipartite graph.}
\label{fig:RP2Graph}
\end{subfigure}

\caption{A triangulation of $\mathbb{R}P^2$ and the associated bipartite graph illustrating Example \ref{projplane}.}
\label{fig:RP2}
\end{figure}

\begin{example}\label{projplane}  Let $X \sim_h \Sigma \mathbb{R}P^2$. 
Consider $Y$ to be the triangulation of $\mathbb{R}P^2$ shown in Figure~\ref{fig:RP2Triangulation}. Then the graph $G$ has $V(G) = \{0, 1, 2, 3, 4, 5\} \sqcup \{
012, 123, 025, 014, 234, 245, 035, 034, 145, 135\},$
\noindent with an edge connecting $i$ ($0\leq i \leq 5$) to $jkl$ whenever $i \not \in \{ j,k, l\}$. The graph $G$ is shown in Figure~\ref{fig:RP2Graph}.
\end{example}

One might hope this method could provide potential counterexamples to Conjectures \ref{conj1} and~\ref{conj2}. Namely, we pick a simplicial complex $Y$ with torsion in its homology. If the corresponding graph $G(Y)$ is a circle graph, it would yield an example of an independence complex of a circle graph that is not homotopy equivalent to a wedge of spheres, since the cohomology of such a wedge of spheres does not contain torsion (recall that $H_{k+1}(\Sigma Y)\cong H_k(Y)$).
However, we will show the graph in Figure  \ref{fig:RP2Graph} is not a circle graph. This follows from Theorem \ref{Sec2MainThm} below which states that Jonsson's construction never yields a circle graph as long as $Y$ satisfies the following mild condition. \\

\noindent
{\bf Vertex Separation Condition:} A simplicial complex $Y$ is said to satisfy the {\it Vertex Separation Condition} if $Y$ contains four facets\footnote{That is, maximal simplices.} $\Delta_i$ ($i= 0, 1, 2, 3$) and three vertices $a_j$ ($j=1,2,3$), such that $a_i \in \Delta_i$, and $a_i \notin \Delta_j$ for $j\neq i$.\\

\begin{figure}[ht]
\centering
\begin{subfigure}{.49\textwidth}
\centering	
\includegraphics[width=0.5\linewidth]{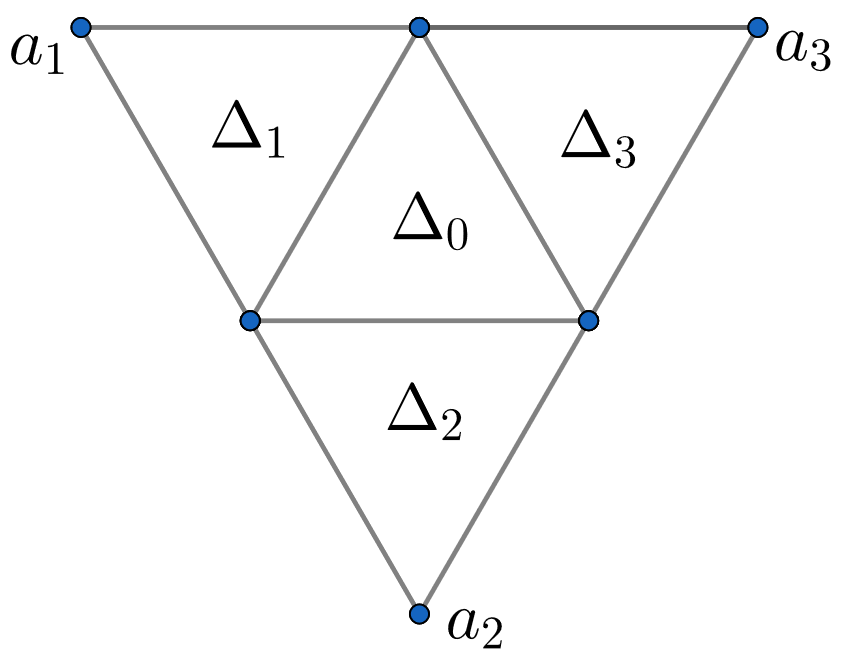}
\caption{Four facets in $Y$.}
\label{fig:4facets}
\end{subfigure}
\begin{subfigure}{.49\textwidth}
\centering
\includegraphics[width=0.5\linewidth]{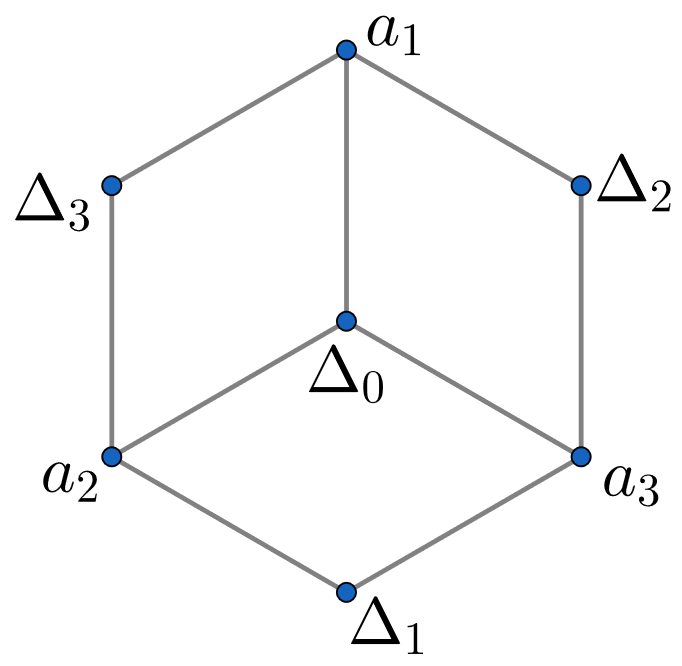}
\caption{Three parallelograms in $G$.}
\label{fig:3parallelograms}
\end{subfigure}
\caption{Four facets and three parallelograms illustrating proof of Theorem \ref{Sec2MainThm}.}
\label{fig:mild_condition}
\end{figure}

This condition is satisfied by the triangulation of $\mathbb{R}P^2$ shown in Figure \ref{fig:RP2Triangulation}, as it contains the pattern in Figure \ref{fig:4facets}.  As we will see in the proof of Theorem \ref{Sec2MainThm}, such a pattern implies that the graph $G$ contains Figure \ref{fig:3parallelograms} as a subgraph.

\begin{theorem} 
\label{Sec2MainThm}
Let $Y$ be a simplicial complex and $G=G(Y)$ be the bipartite graph constructed as above. If $Y$ satisfy the Vertex Separation Condition, then $G$ is not a circle graph. 
\end{theorem}
\begin{proof}
 Consider the seven vertices in $G$ corresponding to $a_j$ \ $(j=1,2,3)$ and $\Delta_i$ \ $(i=0,1,2,3)$. For simplicity, we continue to denote these vertices by $a_j$ and $\Delta_i$.\\
Assume $G$ is a circle graph. Without loss of generality, we assume that the chords associated with the $\Delta_i$'s lie outside the circle, while the chords corresponding to the $a_j$'s lie inside the circle. These chords will be denoted by $\overline{\Delta}_i$ and $\overline{a}_j$ respectively.\\
Since $a_j \notin \Delta_0$ in $Y$ for all $j$, the vertices $a_1, a_2, a_3$ are each adjacent to the vertex $\Delta_0$ in $G$. Therefore, the corresponding chords $\overline{a}_1, \overline{a}_2, \overline{a}_3$ must all intersect the chord $\overline{\Delta}_0$. This implies that the three chords $\overline{a}_1, \overline{a}_2, \overline{a}_3$ must be parallel to each other (see Figure \ref{fig:ChordDiagram}). 
 \begin{figure}[ht]
\centering	
\includegraphics[width=0.4\linewidth]{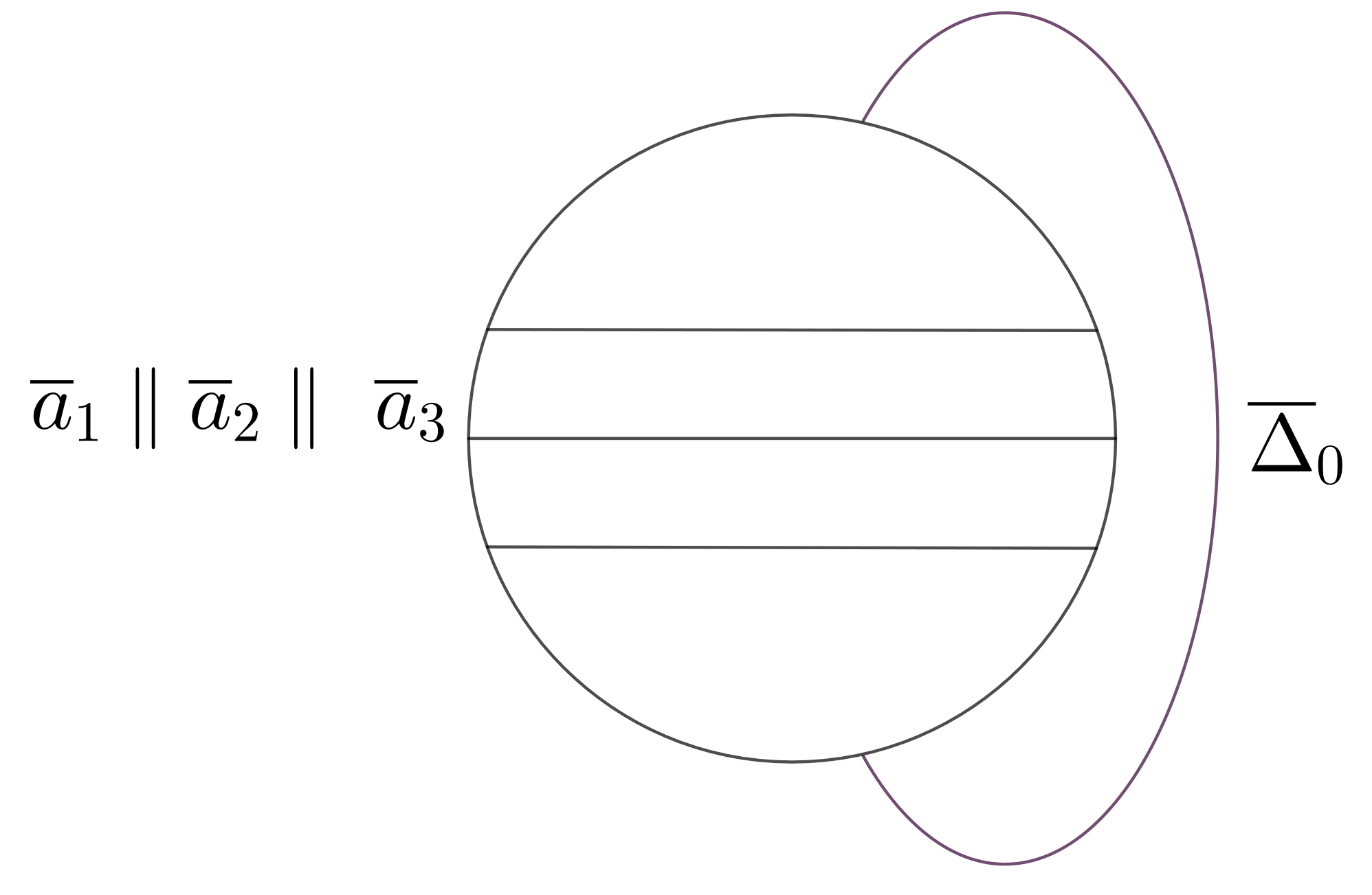}
\caption{Chord Diagram.}
\label{fig:ChordDiagram}
\end{figure}

There are six possible orderings for these chords, corresponding to the six permutations of  $\overline{a}_1, \overline{a}_2$  and $\overline{a}_3$. Next, consider $\Delta_1$ in $Y$, which contains $a_1$ but not $a_2$ or $a_3$. Thus, the corresponding vertex $\Delta_1$ in $G$ is adjacent to the vertices $a_2$ and $a_3$, but not to $a_1$. Therefore, the corresponding chord $\overline{\Delta}_1$ must intersect the chords $\overline{a}_2$ and $\overline{a}_3$, but not $\overline{a}_1$. This implies that the chords $\overline{a}_2$ and $\overline{a}_3$ must be positioned together, meaning that $\overline{a}_1$ cannot be placed between them.

By the same argument, considering $\Delta_2$ (respectively, $\Delta_3$), we conclude that the chord $\overline{a}_2$ (respectively, $\overline{a}_3$) cannot be placed in the middle. Therefore, none of the three chords $\overline{a}_1, \overline{a}_2,$ or $\overline{a}_3$ can be placed in the middle, which leads to a contradiction. This concludes that the graph $G$ cannot be a circle graph.    
\end{proof}

\begin{remark}
\label{RemarkonMildCondisions}
(a) The Vertex Separation Condition  is easily met by many triangulations. For example, if there are four disjoint facets $\Delta_i$ ($i=0,1,2,3)$, then we can simply pick any vertex $a_i$ on $\Delta_i$ ($i=1,2,3$). Or if there is a pattern like in Figure \ref{fig:4facets}, then the condition is satisfied.  In the case of the above triangulation of $\mathbb{R}P^2$, it contains a pattern as Figure \ref{fig:4facets}, thus the graph $G$ cannot be a circle graph.\\
(b) On the other hand, if the condition in Theorem \ref{Sec2MainThm} is not satisfied, then in principle the graph $G$ could be a circle graph. For example, if the simplicial complex $Y$ is the tetrahedron, the graph $G(Y)$ is a circle graph.  This does not provide a counterexample to Conjecture \ref{conj1} because the space $Y$ is homotopy equivalent to $S^2$.
\end{remark}

\section{Examples of Khovanov A-adequate links}\label{Examples}
In this section, we delve into the concept of Khovanov $A$-adequacy by providing and analyzing new examples of families of Khovanov $A$-adequate links. By examining these examples, we aim to enhance our understanding of the underlying structures and properties that characterize Khovanov $A$-adequacy.

\subsection{A family of Khovanov A-adequate twisted torus links}

Consider the braid word $\beta(6,3)=(\sigma_5\sigma_4\sigma_3\sigma_2\sigma_1)^3$. From $\beta(6,3)$ we construct the braid diagram shown in Figure~\ref{fig:TwistedT(3,6)Example} corresponding to the braid word $w=\beta(6,3)\sigma_2\sigma_4(\sigma^{-1}_1\sigma^{-1}_3\sigma^{-1}_5)\sigma_2\sigma_4$. Our first example is the diagram $D$ obtained by closing $w$, whose associated $s_AD$ smoothed diagram is shown in Figure \ref{fig:TwistedT(3,6)Example_SAstate}. 

\begin{figure}[ht]
\centering
\begin{subfigure}{.49\textwidth}
\centering		\includegraphics[width=0.60\linewidth]{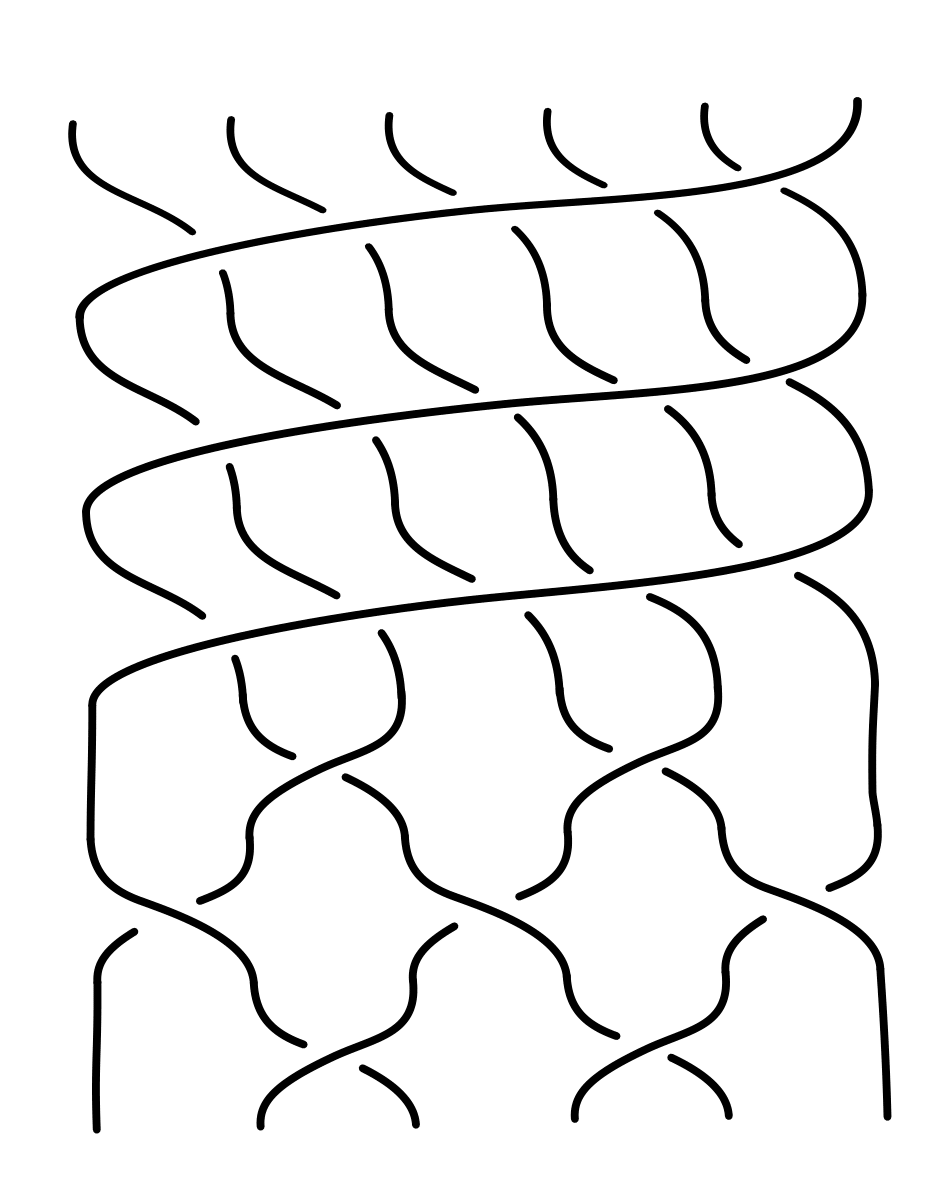}
\caption{A diagram of  $\beta(6,3)\sigma_2\sigma_4(\sigma^{-1}_1\sigma^{-1}_3\sigma^{-1}_5)\sigma_2\sigma_4$.}
\label{fig:TwistedT(3,6)Example}
\end{subfigure}
\begin{subfigure}{.49\textwidth}
\centering
\includegraphics[width=0.60\linewidth]{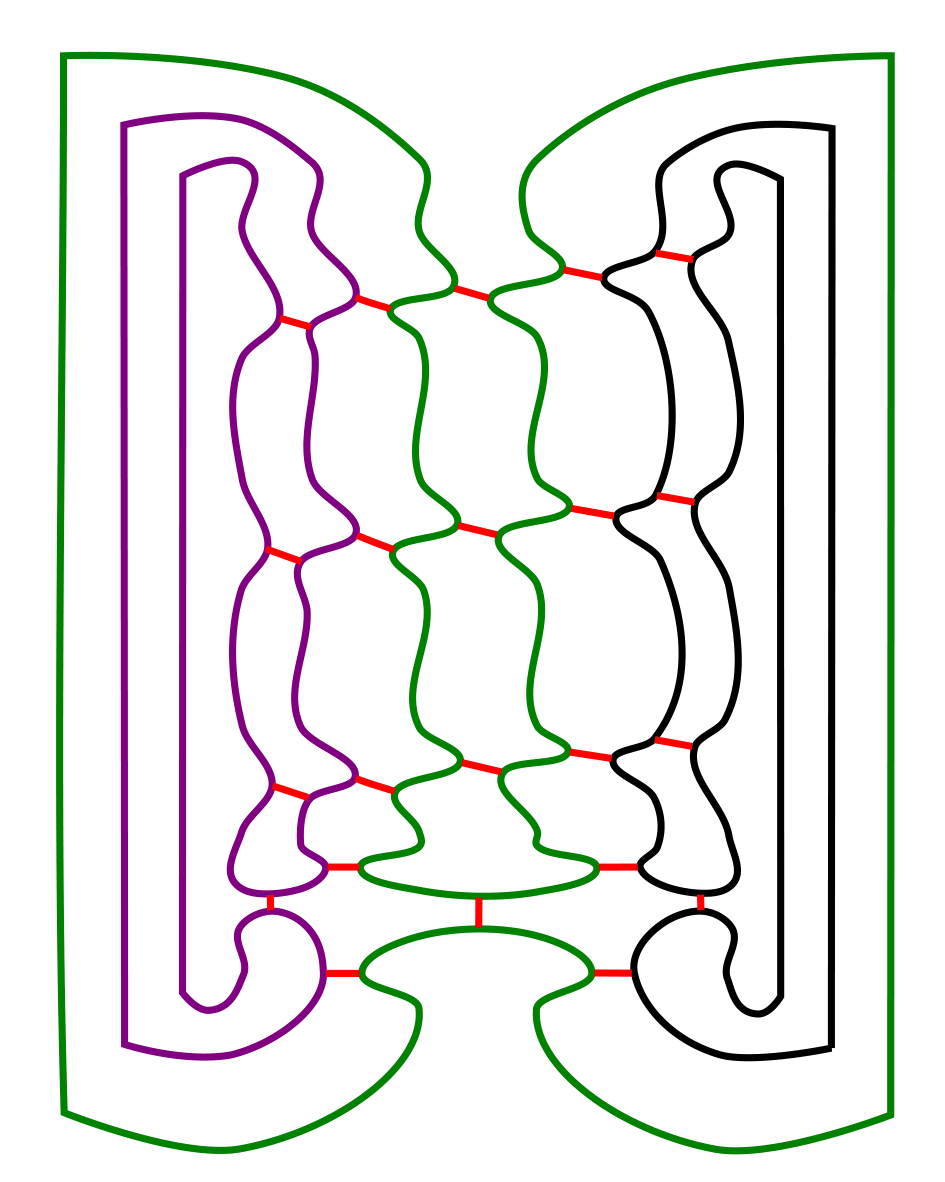}
\caption{The $s_AD$ state associated to its closure.}
\label{fig:TwistedT(3,6)Example_SAstate}
\end{subfigure}
\caption{A braid-like diagram and the associated $s_A$ state of its closure.}
\label{fig:TwistedExample}
\end{figure}

The Lando graph $G_D$ associated to $D$ consists of the disjoint union of three stars, each of them with $3$ rays, as shown in Figure \ref{GraphIndepT(6,3)}. The independence complex associated to each of those stars is homotopy equivalent to $S^0$, and therefore $I(D) = I(G_D)) \simeq S^{0} \ast S^{0} \ast  S^{0} \simeq S^{2}$, which is not contractible. Therefore, the diagram $D$ is Khovanov $A$-adequate, and so is the associated link.

\begin{figure}[h]
    \centering
\includegraphics[width=0.85\linewidth]{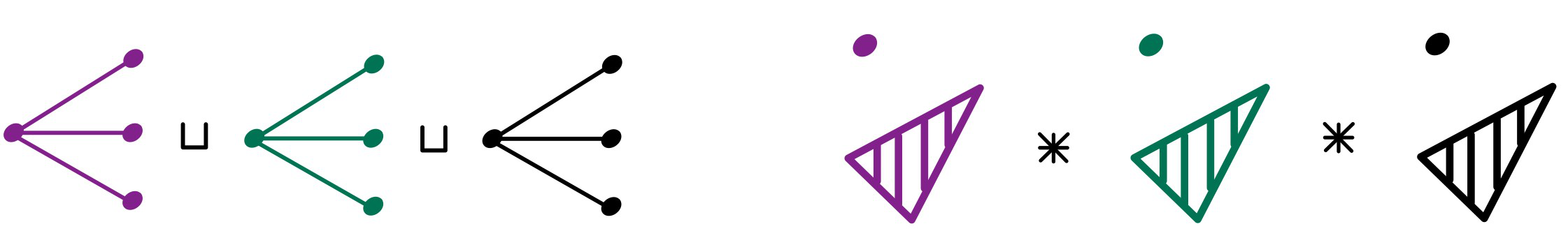}
\caption{Lando graph associated to diagram $D$ (left) and its associated independence complex (right).}
\label{GraphIndepT(6,3)}
\end{figure}


We next generalize the previous example. Given two positive integers $m,n$, with $m \geq 6$, write $\beta(m,n)$ for the braid word whose closure is the standard diagram of the torus link $T(m,n)$, that is, $\beta(m,n) = (\sigma_{m-1} \cdots \sigma_2 \sigma_1)^n$.  

\begin{figure}[h]
    \centering
\includegraphics[width=0.5\linewidth]{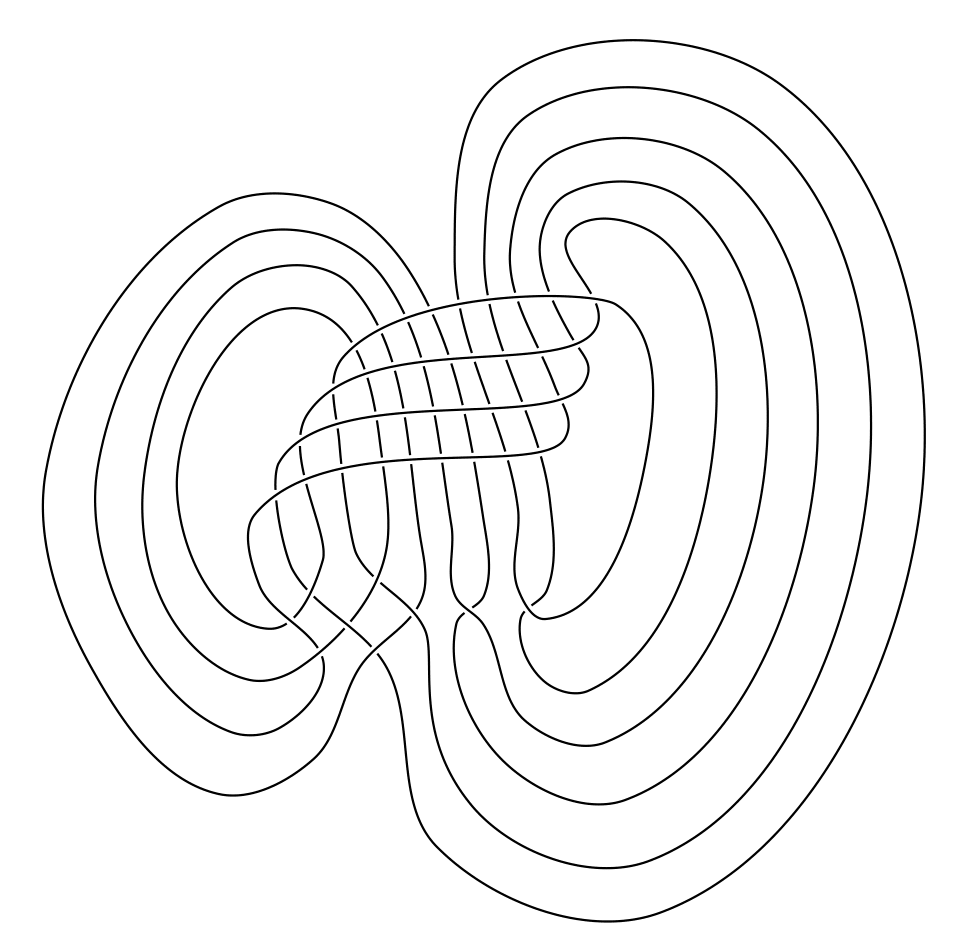}
\caption{Diagram $D_{10,4}$.}
    \label{beta(9,4)}
\end{figure}

\begin{proposition}
Given two positive integers $m,n$, with $m\geq 6$, consider the diagram $D_{m,n}$ obtained as the closure of the braid word $w_{m,n}$ defined as  $$w_{m,n} = \left\{ \begin{array}{lll}   \beta(m,n)\sigma_2\sigma_4(\sigma^{-1}_1\sigma^{-1}_3 \cdots \sigma^{-1}_{m-1})\sigma_2\sigma_4 & & \mbox{ if } m=2k,  \\ \\ \beta(m,n)\sigma_2\sigma_4(\sigma^{-1}_1\sigma^{-1}_3 \cdots \sigma^{-1}_{m-2})\sigma_2\sigma_4 & & \mbox{ if } m = 2k+1. \end{array} \right. $$

Then $I(D_{m,n}) \sim_h S^{k-1}$, and therefore its extreme Khovanov homology is isomorphic to $\mathbb{Z}$.
\end{proposition}

\begin{proof}
If $m$ is even, $s_AD_{m,n}$ consists of $k$ circles, and for each circle $c$ there are four chords having both endpoints on $c$ in a similar pattern as that in Figure \ref{fig:TwistedT(3,6)Example_SAstate}. Therefore, $I(D_{m,n})$ is homotopy equivalent to the join of $k$ copies of $S^0$, and the result holds. 

The case when $m$ is odd is analogous, since $s_AD_{m,n}$ consists of $k+1$ circles but one of the circles is not contributing to the Lando graph $G_{D_{m,n}}$, since there are no chords having both endpoints on it. 
\end{proof}

\subsection{A family of Khovanov $A$-adequate negative braid links}

Our next family of examples consists of the closure of certain negative braid diagrams. It is known that negative (resp. positive) diagrams are $B$-adequate and therefore Khovanov $B$-adequate (resp. $A$-adequate and therefore $A$-adequate). However, there is no reason to expect an arbitrary negative (resp. positive) diagram to be Khovanov $A$-adequate (resp. $B$-adequate).

\ 

In this section we show diagrams that are Khovanov adequate (i.e., they are simultaneously Khovanov $A$-adequate and Khovanov $B$-adequate). In addition, these diagrams are $B$-adequate, but none are $A$-adequate. Therefore, these examples show that the concept of Khovanov adequacy properly extends the classical definition of adequacy at the level of link diagrams.  

\ 

\begin{figure}[ht]
\centering
\begin{subfigure}{.49\textwidth}
\centering		
\includegraphics[height=7cm]{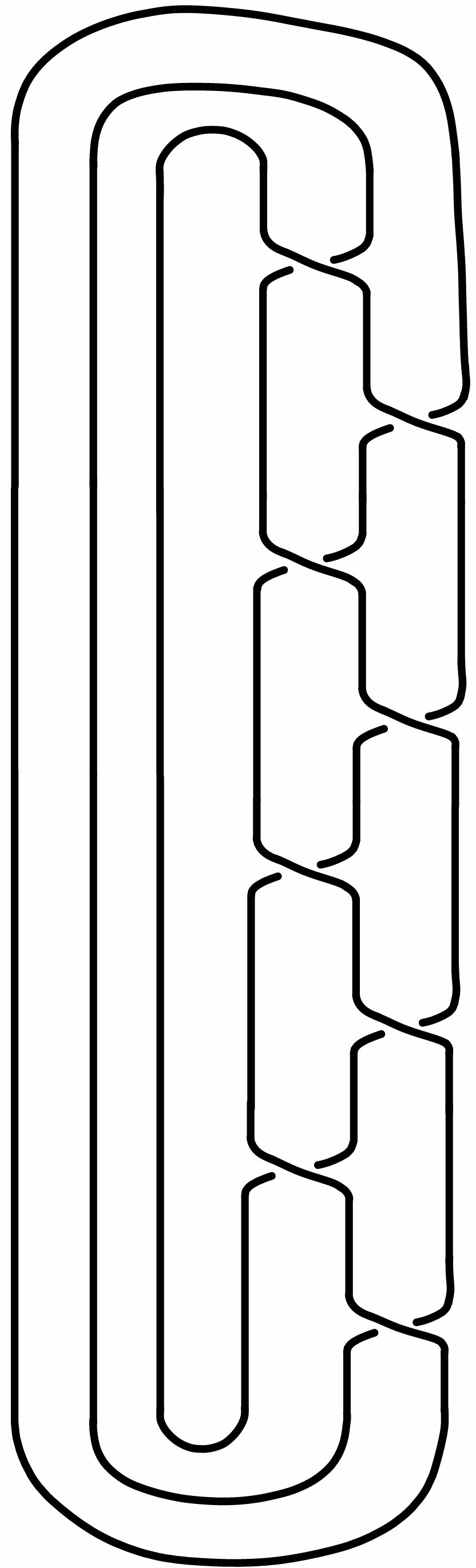}
\caption{Standard diagram of $T(3,-4)$ representing $8_{19}$.}
\label{fig:819 diagram}
\end{subfigure}
\begin{subfigure}{.49\textwidth}
\centering
\includegraphics[height=7cm]{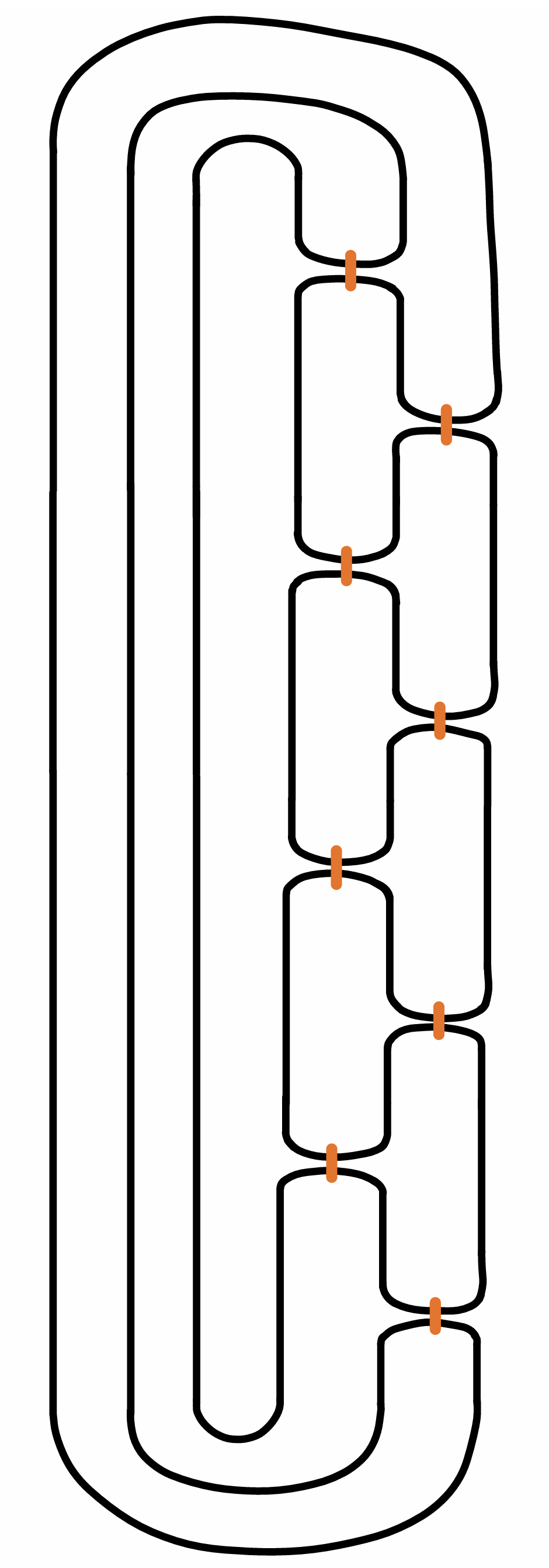}
\caption{The associated smoothed state $s_AD$.}
\label{fig:819 A-smoothed}
\end{subfigure}
\caption{A diagram $D$ of knot $8_{19}$ and $s_AD$.}
\label{fig:819}
\end{figure}

Our starting example is a diagram $D$ of the knot $8_{19}$ drawn as the standard diagram of the torus knot $T(3, -4)$, as shown in Figure \ref{fig:819}. Its Lando graph $G_D$ is the $8$-gon $C_{8}$, and thus by Proposition~\ref{cycle pattern lemma} we get that $I(D)) \sim_h S^2$. 

By increasing $|r|$ we recover the following result from \cite[Section 7.3]{PS1}: 
\begin{proposition} \label{prop C2r}
Let $D_r$ be the standard diagram of the torus link $T(3,-r)$. Then, 
$$I(D_r) \simeq_h \left\{ \begin{array}{lll}   S^{k-1} & & \mbox{ if } 2r=3k\pm 1, \\ \\ S^{k-1} \vee S^{k-1} & & \mbox{ if } 2r = 3k. \end{array} \right. $$
\end{proposition}

In \cite{PS1} Przytycki and Silvero computed (real) extreme Khovanov homology of torus links of type $T(3,-r)$. However, we include these examples as the starting point to construct an infinite family of braid diagrams with a larger number of strands, whose Lando graphs are identical to those considered in this section. 



\subsubsection{The Family $F(s,r)$}

\begin{figure}[ht]
\centering
\begin{subfigure}{.49\textwidth}
\centering		
\includegraphics[height=7cm]{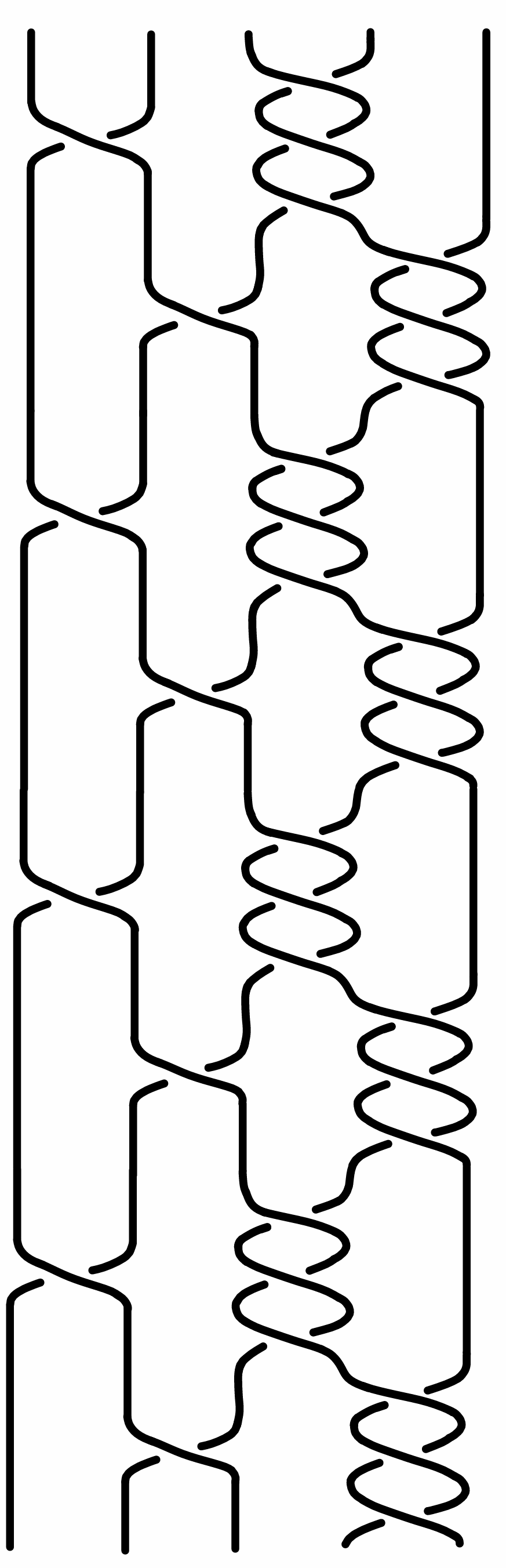}
\caption{Knot diagram $D=F(5, 4)$. }
\label{fig:example 2 diagram}
\end{subfigure}
\begin{subfigure}{.49\textwidth}
\centering
\includegraphics[height=7cm]{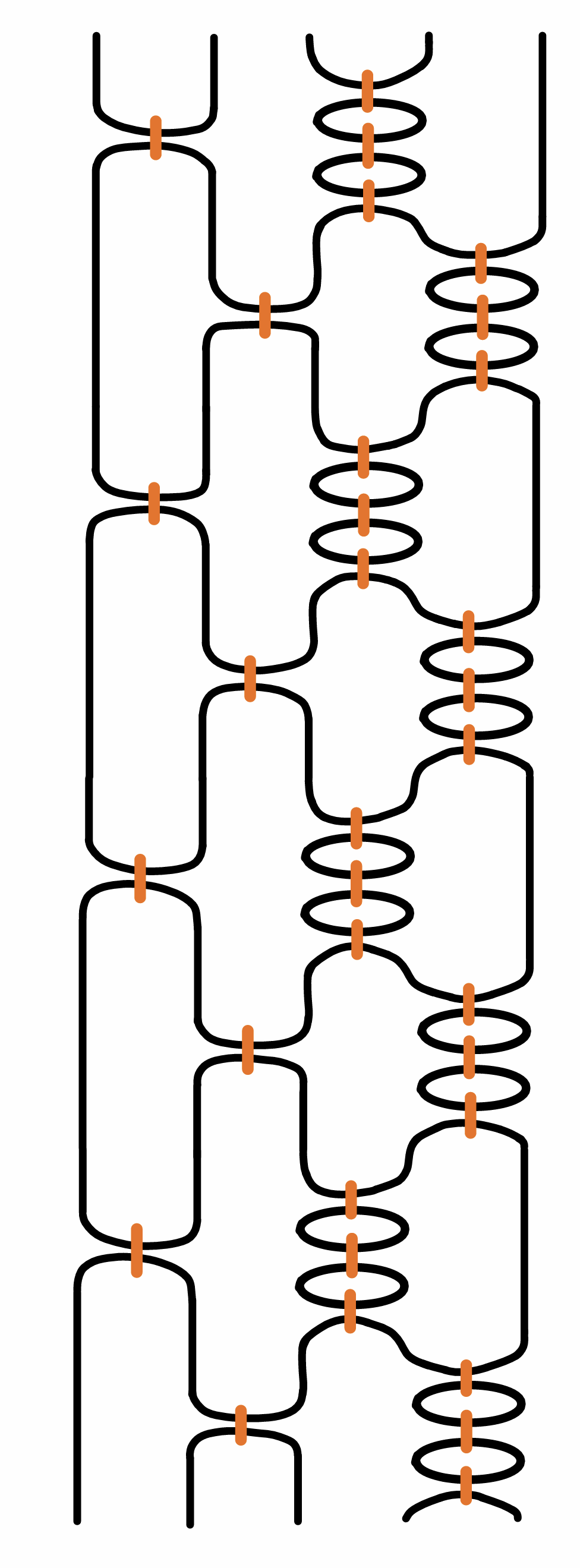}
\caption{The associated smoothed state $s_AD$.}
\label{fig:example 2 diagram A-smoothed}
\end{subfigure}
\caption{Diagram $D= F(5,4)$ and the associated smoothed diagram $s_AD$.}
\label{fig:example 2}
\end{figure}

Given a positive integer $r$ and an odd positive integer $s \geq 5$, we write $F(s,r)$ for the link diagram consisting of the closure of the braid diagram associated to the word $(\sigma_1^{-1} \sigma_3^{-3} \dots \sigma_{s-2}^{-3} \sigma_2^{-1} \sigma_4^{-3} \dots \sigma_{s-1}^{-3})^r$. Diagram $F(5,4)$ is depicted in Figure \ref{fig:example 2 diagram}. 


If we compare Figures \ref{fig:819 A-smoothed} and \ref{fig:example 2 diagram A-smoothed} we see that the arrangement of chords with both endpoints in the same circle associated to $F(5,4)$ and the standard diagram of $T(3,-4)$ is the same, and therefore the Lando graph associated to both diagrams coincide. This observation can be generalized as follows:





\begin{proposition}
In the previous setting, the Lando graph associated to $F(s,r)$ is the $2r$-gon, and therefore 

$$I(F(s,r)) \sim_h  \begin{cases} 
    S^{k-1} & 2r = 3k \pm 1,\\
    S^{k-1}\vee S^{k-1} & 2r = 3k.
    \end{cases}$$
\end{proposition}   

Since the independence complexes associated to diagrams $F(s,r)$ are not contractible, these are examples of Khovanov $A$-adequate diagrams.

\subsection{Modified cable link diagrams}

Roughly speaking, a cable link of a knot $K$ is a satellite knot with pattern a torus knot and companion $K$. Given an $A$-adequate knot diagram $D$, any associated cable diagram with no twistings (i.e., built as a set of parallel copies of $D$) is still $A$-adequate. However, we observe that the diagram obtained from such a cable link after introducing a local Hopf tangle is not necessarily $A$-adequate, and it can, in fact, be Khovanov $A$-adequate.  

\begin{figure}[ht]
\centering
\begin{subfigure}{.49\textwidth}
\centering		\includegraphics[width=0.7\linewidth]{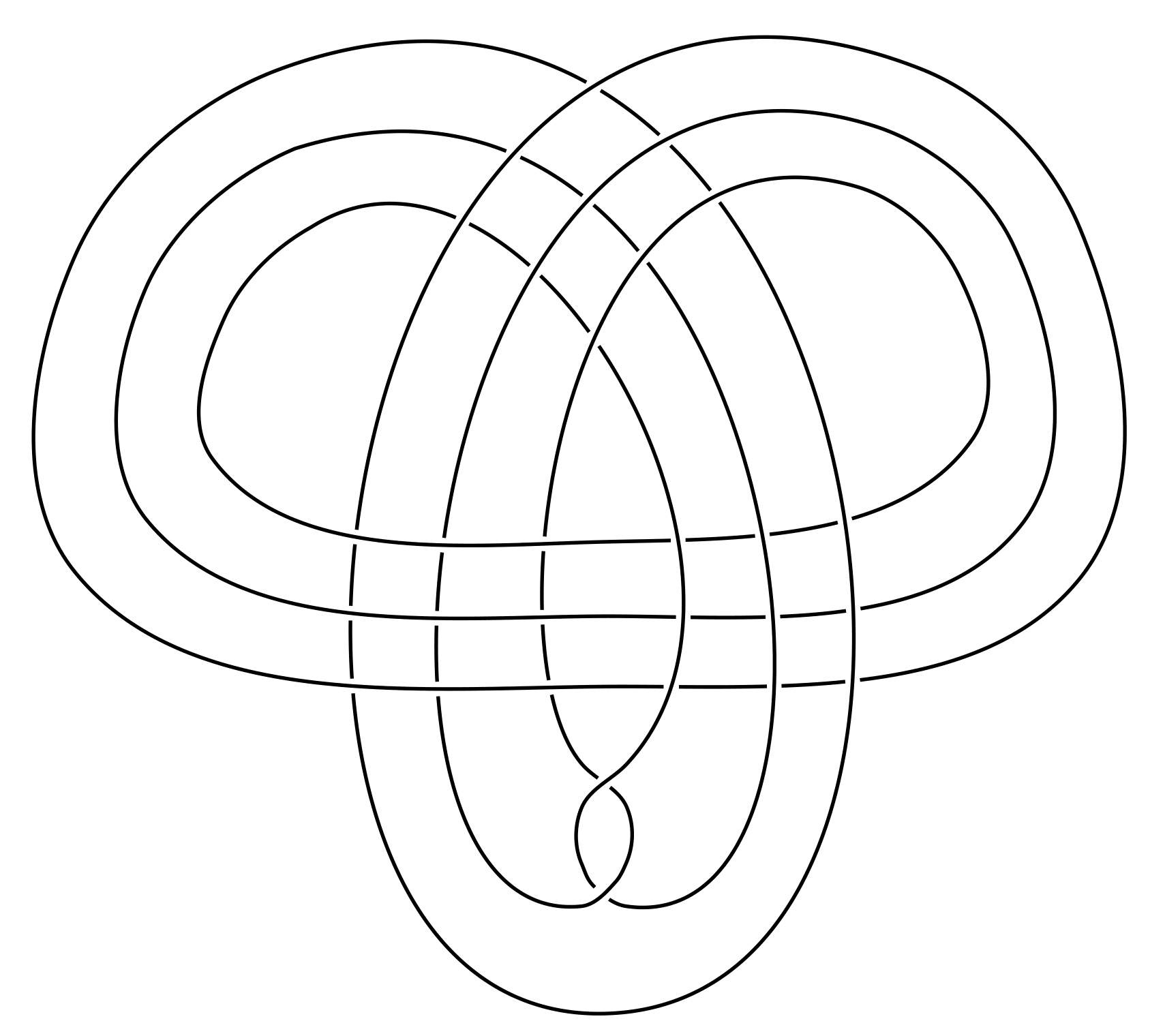}
\caption{A $3$-cable of the trefoil after adding a Hopf $2$-tangle.}
\label{Figura1.png}
\end{subfigure}
\begin{subfigure}{.49\textwidth}
\centering
\includegraphics[width=0.7\linewidth]{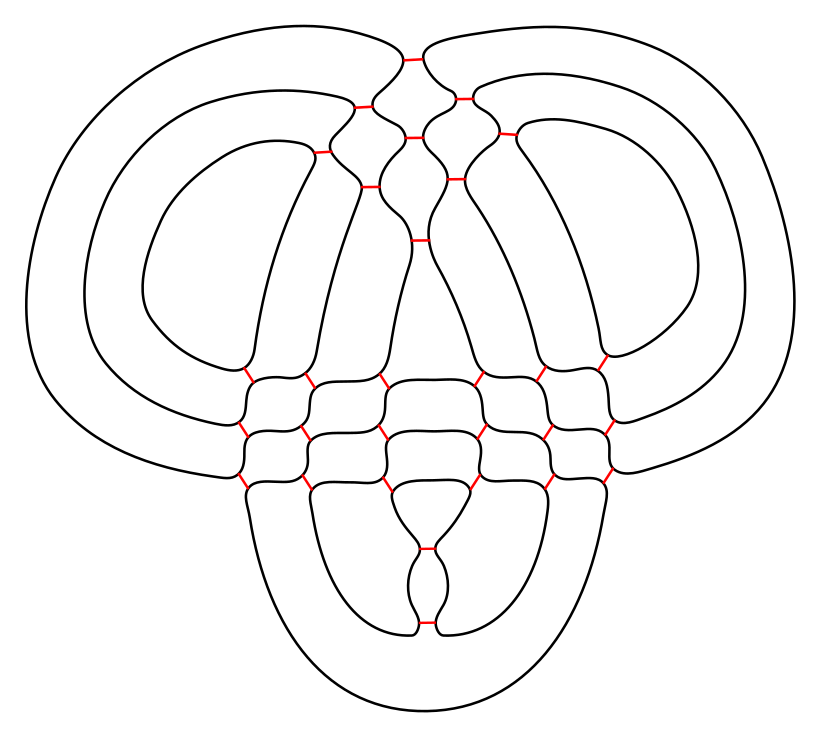}
\caption{The associated smoothed state $s_AD$.}
\label{Figura1_A-smoothing(1).png}
\end{subfigure}
\caption{An example of a Khovanov $A$-adequate diagram which is not $A$-adequate.}
\label{figsAli}
\end{figure}

We show such an example: consider the diagram $D$ shown in Figure \ref{Figura1.png}, which consists of a $3$-cable of the trefoil knot together with a Hopf tangle. It follows from Figure \ref{Figura1_A-smoothing(1).png} that  the Lando graph $G_D$ is the $4$-gon, and therefore $I(D)\sim_h S^0$. Hence, $D$ is a Khovanov $A$-adequate diagram which is not $A$-adequate.

\

\section*{Acknowledgements}
This paper arises as the result of the work carried out by the authors during their participation in the program ``Low Dimensional Topology: Invariants of Links, Homology Theories, and Complexity'' held at the Matrix Research Institute (The University of Melbourne) in June 2024. We are grateful to the center for their hospitality. LB was partially supported by NSF Grant DMS-2038103. GMV was supported by the Matrix-Simons travel grant and acknowledges the support of the National Science Foundation through Grant DMS-2212736. YR was supported by the Matrix-Simons travel grant. MS was partially supported by the Spanish Research Grant PID2020-117971GB-C21 funded by MCIN/AEI/10.13039/501100011033.

\end{document}